\documentclass{commat}

\usepackage{enumitem}

\title{%
    Positive solutions for a Schr\"odinger-Bopp-Podolsky system
    }

\author{%
    Bruno Mascaro and Gaetano Siciliano
    }

\affiliation{
    \address{Bruno Mascaro -
    Departamento de Matem\'atica, 
	Instituto de Matemática e Estat\'istica, 
	Universidade de Sa\~o Paulo,
	Rua do Mat\~ao 1010,  05508-090 S\~ao Paulo, SP, Brazil
        }
    \email{%
    bmascaro@ime.usp.br
    }
    \address{Gaetano Siciliano -
     Departamento de Matem\'atica, 
		Instituto de Matem\'atica e Estat\'istica, 
		Universidade de S\~ao Paulo,
		Rua do Mat\~ao 1010,  05508-090 S\~ao Paulo, SP, Brazil
        }
    \email{%
    sicilian@ime.usp.br
    }
    }

\abstract{%
    We consider the following Schr\"odinger-Bopp-Podolsky system in $\mathbb R^{3}$
$$\left\{ \begin{array}[c]{ll}
-\varepsilon^{2} \Delta u + V(x)u + \phi u  =  f(u) \medskip\\ 
-\varepsilon^{2} \Delta \phi + \varepsilon^{4} \Delta^{2}\phi  =  4\pi u^{2},
\end{array}
\right.$$
where $\varepsilon > 0$ with $ V:\mathbb{R}^{3} \rightarrow \mathbb{R}, f:\mathbb{R} \rightarrow \mathbb{R}$ 
satisfy suitable assumptions.
By using variational methods, we prove that the number of positive solutions is estimated below by the Ljusternick-Schnirelmann category of $M$, the set of minima of the potential~$V$.
    }

\keywords{%
   Variational methods, Ljusternick-Schnirelmann category, multiplicity of solutions.
    }

\msc{%
    35A15, 35J58, 35Q55, 35Q61
    }

\VOLUME{31}
\YEAR{2023}
\NUMBER{1}
\firstpage{237}
\DOI{https://doi.org/10.46298/cm.10363}

\begin{paper}

\section{Introduction}

In this paper, by using the ideas
developed originally in 
\cite{11}, \cite{BC91}, \cite{z4} we show existence and multiplicity of positive solutions for the following
problem in $\mathbb R^{3}$
\begin{equation}\label{eq:Pe}\tag{$P_{\varepsilon}$}
\left\{
 \begin{array}[c]{ll}
-\Delta u +V(\varepsilon x)u + \phi u  = f(u) \medskip \\
-\Delta \phi  + \Delta^{2}\phi  = 4\pi u^{2}
\end{array}
\right.
\end{equation}
whenever $\varepsilon>0$ is a small parameter and  $f$ and $V$
satisfy the   assumptions given below.

In the mathematical literature such  a problem was introduced recently in
\cite{3} and  describes the stationary states of the Schr\"odinger equation
in the generalised electrodynamics developed by  Bopp and Podolsky.
Roughly speaking, the system appears when one searches for stationary solutions,
namely solutions of type $\psi(x,t)= u(x)e^{it}$, of the Schr\"odinger equation of a moving charged particle
which interacts with its own purely electrostatic field, in the case in which the generalised electromagnetic theory of 
Bopp-Podolsky is considered. The reason to prefer the Bopp-Podolsky theory to the classical and more studied Maxwell
theory, is that in the first case the energy associated to a charged particle is finite.
In fact in the Bopp-Podolsky generalised electrodynamic the equation 
of the 
electrostatic field generated by a charge particle (let us say at rest in the origin) is
$$-\Delta \phi +\Delta^{2} \phi = \delta$$
and the fundamental solution $\Phi(x) = \frac{1- e^{-1/|x|}}{|x|}$
has finite energy, $\int |\nabla \Phi|^{2} + \int |\Delta \Phi|^{2} <+\infty$.
On the other hand  in the Maxwell theory the equation of the electrostatic field
is
$$-\Delta \phi = \delta$$
and the fundamental solution $\Phi(x) = \frac{1}{|x|}$ satisfies $\int |\nabla \Phi|^{2}=+\infty$,
giving rise to the so-called {\sl infinity problem}. 
 Then the equations in system \eqref{eq:Pe} have to be interpreted as 
 \begin{itemize}
 \item a Schr\"odinger type equation
 (the first one) in presence of a fixed external potential $V$ and an ``internal'' potential $\phi$, and
 \item an  equation (the second one) which says that the potential $\phi$ has as source the same wave function, being $|\psi| = u^{2}$,  justifying the term ``internal''.
 \end{itemize}

For the mathematical derivation of such a system and some related results
we refer the reader to the recent papers \cite{3}, \cite{GK}, \cite{GD}, \cite{GG}
where the problem is studied under various conditions or even in a bounded domain.

In this paper we assume that  $V$ and the nonlinearity $f$ satisfy

 \begin{enumerate}[label=(V\arabic*),ref=V\arabic*,start=1]
 \item\label{V} $V:\mathbb{R}^{3} \rightarrow \mathbb{R}$ is a continuous function 
 such that
 $$0< \min_{\mathbb{R}^{3}} V:=V_{0}<V_{\infty}:=\liminf_{|x| \rightarrow +\infty} V \in (V_{0},+\infty],$$
with $M=\{ x\in \mathbb{R}^{3}:V(x) = V_{0} \}$ smooth and bounded,
 \end{enumerate}

\begin{enumerate}[label=(f\arabic*),ref=f\arabic*,start=1]
\item\label{f1} $f:\mathbb{R} \rightarrow \mathbb{R}$ is a function of class $C^{1}$ and $f(t) = 0$ for $t=0$,
\item\label{f2} $\lim_{t\rightarrow 0} \frac{f(t)}{t} = 0$,
\item\label{f3} there exists $q_{0} \in (3, 2^{*}-1)$ such that 
$\lim_{t\rightarrow +\infty} \frac{f(t)}{t^{q_{0}}} = 0$, where $2^{*}  = 6$,
\item\label{f4} there exists $K>4$ such that $0<KF(t):= K\int_{0}^{t}f(\tau)d\tau \leq tf(t)$ for all $t>0$,
\item\label{f5} the function $t\mapsto \frac{f(t)}{t^{3}}$ is strictly increasing in $(0,+\infty)$.
 \end{enumerate}

The assumptions on the nonlinearity $f$ are standards in order to work with variational methods,
use the Nehari manifold and deal with the Palais-Smale condition.

The assumption \eqref{V} will be fundamental in order to estimate the number of positive solutions and also to recover some compactness.

The main result of this work is:
\begin{theorem}\label{th:main}
Under the above assumptions \eqref{V}, \eqref{f1}-\eqref{f5}, there exists an $\varepsilon^{*}>0$ such that for every $\varepsilon \in (0, \varepsilon^{*}]$,  problem \eqref{eq:Pe} possesses at least $cat M$ positive solutions. Moreover, if $cat M >1$, then (for a suitably small $\varepsilon$) there exist at least $cat M + 1$ positive solutions.
\end{theorem}

In particular among these solutions there is the ground state, namely the solution with minimal energy;
this will be evident by the proof.
Here $ \text{cat} M:=\text{cat}_{M}M $ is the Ljusternik-Schnirelmann category and by positive solutions 
we mean a pair $(u,\phi)$ with $u$ positive, since $\phi$
will be automatically positive.

We point out that the assumption on the potential $V$ is not too restrictive since it is satisfied by
an interesting class which appear in physical models, such as the confining potentials.
There is then an interesting relation between the topology of the set of minima of $V$
and the  number of solutions.

For example, for a potential of type
$$
V(x) = 
\begin{cases}
1 & \mbox{ if } |x| \leq1 \\
|x|^{2}&\mbox{ otherwise } 
\end{cases}
$$
the theorem states the existence of (at least) one solution for small $\varepsilon$,
 being $M$ the unit ball $ \{ x\in \mathbb R^{3} : |x|\leq1\}$ and $cat M=1$. On the other hand with the
 following double-well potential
$$
V(x) = 
\begin{cases}
1 & \mbox{ if } \ \pi/2 \leq |x|\leq 3\pi/2\ \   \text{ and } \ \ \ 5\pi/2 \leq |x|\leq 7\pi/2 \\
2+ \cos|x| & \mbox{ if }\  0\leq |x|\leq \pi/2\ \   \text{ and } \ \ \ 3\pi/2 \leq |x|\leq 5\pi/2\\
|x|^{2}+1-49\pi^{2}/4&\mbox{ if } \ |x|\geq7\pi/2
\end{cases}
$$
the theorem states that there are at least three solutions for small $\varepsilon$, since
 $M$ is  the union of the annuli $ \{ x\in \mathbb R^{3} : \pi/2\leq  |x|\leq 3\pi/2\}$ and
$ \{ x\in \mathbb R^{3} : 5\pi/2\leq  |x|\leq 7\pi/2\} $ 
and $cat M=2$.

As a matter of notations, all the integrals, unless otherwise specified, are understood on $\mathbb R^{3}$
with the Lebesque measure. We denote with $\|\cdot\|_{p}$ the usual $L^{p}$ norm.
Finally $o_{n}(1)$ denotes a vanishing sequence and we use the letter $C$ to denote a 
positive constant whose value does not matter and can vary from line to line.

\section{Preliminaries and Variational Setting}
Let us start by  recalling some results that will be useful for our work. For more details 
 see \cite{3}.

Let $\mathcal{D}$ be the completion of $C^{\infty}_{0}(\mathbb{R}^{3})$ with respect to the norm $\|\cdot\|_{\mathcal{D}}$ induced by $$<\phi ,\psi>_{\mathcal{D}} = \int{\nabla \phi \nabla \psi} + \int{\Delta \phi \Delta \psi}.$$
The space $\mathcal{D}$ is an Hilbert space continuously embedded into $\mathcal{D}^{1,2}(\mathbb{R}^{3})$ and consequently in $L^{6}(\mathbb{R}^{3}).$
Moreover this space is embedded also into $L^{\infty}(\mathbb R^{3})$.

The following lemmas are used to justify a ``reduction method''
in order to  deal with just one equation.


\begin{lemma}
The space $C^{\infty}_{0}(\mathbb{R}^{3})$ is dense in $$A = \{ \psi \in \mathcal{D}^{1,2}(\mathbb{R}^{3}) : \Delta \psi \in L^{2}(\mathbb{R}^{3})\}$$
normed by $\sqrt{<\phi, \phi>_{\mathcal{D}}}$ and, therefore, $\mathcal{D}=A$.
\end{lemma}

For every fixed $u \in H^{1}(\mathbb{R}^{3})$, the Riesz theorem implies that there exists a unique solution $\phi_{u} \in \mathcal{D}$, for the second equation in \eqref{eq:Pe}. Such a solution is given by $\phi_{u} = K*u^{2},$
where $$K(x) =  \frac{1-e^{-|x|}}{|x|}.$$
%
%
%
The solution $\phi_{u}$ has the following properties.
\begin{lemma}
For every $u \in H^{1}(\mathbb{R}^{3})$ we have:
\begin{itemize}
\item[i)] $\forall y \in \mathbb{R}^{3}, \phi_{u(.+y)} = \phi_{u}(.+y)$;
\item[ii)] $\phi_{u} \geq 0$;
\item[iii)]  $\forall s \in (3, +\infty], \phi_{u} \in L^{s}(\mathbb{R}^{3})\cap C_{0}(\mathbb{R}^{3})$;
\item[iv)]$\forall s \in ({3}/{2}, +\infty], \nabla \phi_{u} = \nabla K*u^{2} \in L^{s}(\mathbb{R}^{3})\cap C_{0}(\mathbb{R}^{3})$;
\item[v)] $\|\phi_{u}\|_{6} \leq C \|u\|^{2}$;
\item[vi)] $\phi_{u}$ is the unique minimizer of the functional $$E(\phi) = \frac{1}{2}\|\nabla \phi\|^{2}_{2} + \frac{1}{2} \|\Delta \phi\|^{2}_{2} - \int{\phi u^{2}}, \qquad \phi \in \mathcal{D}.$$
Moreover,
\item[vii)] if $v_{n} \rightharpoonup v$ in $H^{1}(\mathbb{R}^{3})$, then $\phi_{v_{n}} \rightharpoonup \phi_{v}$ in $\mathcal{D}$.
\end{itemize}
\end{lemma}


It is easy to see that after a change of variable our problem can be written as
\begin{equation}\label{eq:Pe*}\tag{$P_{\varepsilon}^{*}$}
\left\{
 \begin{array}[c]{ll}
-\Delta u +V(\varepsilon x)u + \phi u  = f(u), \medskip \\
-\Delta \phi  + \Delta^{2}\phi  = 4\pi u^{2}.
\end{array}
\right.
\end{equation}
Hence the critical points of the functional
$$\mathcal I_{\varepsilon}(u,\phi) = \frac{1}{2}\|u\|^{2}_{2} + \frac{1}{2} \int{V(\varepsilon x)u^{2}} +\frac{1}{2}\int{\phi u^{2}} - \frac{1}{16\pi}\|\nabla\phi\|^{2}_{2} - \frac{1}{16\pi}\|\Delta \phi\|^{2}_{2} - \int{F(u)}$$
in $H^{1}(\mathbb{R}^{3})\times \mathcal{D}$ are easily seen to be weak solutions of \eqref{eq:Pe*}; indeed
such a critical point  $(u,\phi) \in H^{1}(\mathbb{R}^{3}) \times \mathcal{D}$ satisfies 
$$0=\partial_{u}\mathcal I_{\varepsilon}(u,\phi)[v] = \int{\nabla u\nabla v}+\int{V(\varepsilon x)uv} + \int{\phi uv}-\int{f(u)u}, \qquad v\in H^{1}(\mathbb{R}^{3}),$$
$$0 = \partial_{\phi}\mathcal I_{\varepsilon}(u,\phi)[\xi] = \frac{1}{2} \int{u^{2}\xi} - \frac{1}{8\pi}\int{\nabla \phi \nabla \xi} - \frac{1}{8\pi}\int{\Delta \phi \Delta \xi}, \qquad \xi \in \mathcal{D}.$$
The next step is  the usual reduction argument in order to deal  with a one variable functional.
Noting that $\partial_{\phi}\mathcal I_{\varepsilon}$ is a $C^{1}$ function and defining $G_{\Phi}$ as the graph of the map $\Phi : u\in H^{1}(\mathbb{R}^{3}) \mapsto \phi_{u} \in \mathcal{D}$, an application of the Implicit Function Theorem gives
$$G_{\Phi} = \{ (u,\phi) \in H^{1}(\mathbb{R}^{3}) \times D: \partial_{\phi} \mathcal I_{\varepsilon}(u,\phi)=0\}, \qquad  \Phi \in C^{1}(H^{1}(\mathbb{R}^{3}),\mathcal{D}).$$
Then
\begin{eqnarray*}
0 = \partial_{\phi}\mathcal I_{\varepsilon}(u,\Phi(u)) = \frac{1}{2} \int{\phi_{u}u^{2}} - \frac{1}{8\pi}\|\nabla \phi_{u}\|^{2}_{2} - \frac{1}{8\pi}\|\Delta \phi_{u}\|^{2}_{2}
\end{eqnarray*}
and substituting
$$
-\frac{1}{4}\int{\phi_{u}u^{2}} = -\frac{1}{16\pi}\|\nabla \phi_{u}\|^{2}_{2} - \frac{1}{16\pi}\|\Delta \phi_{u}\|^{2}_{2}
$$
in the expression of $\mathcal I_{\varepsilon}$ we obtain  the  functional
$$I_{\varepsilon}(u) := \mathcal I_{\varepsilon}(u,\Phi(u)) = \frac{1}{2}\|\nabla u\|^{2}_{2} + \frac{1}{2}\int{V(\varepsilon x)u^{2}} + \frac{1}{4}\int{\phi_{u}u^{2}} - \int{F(u)}.$$
This functional is of class $C^{1}$ in $H^{1}(\mathbb{R}^{3})$ and, for all $u,v \in H^{1}(\mathbb{R}^{3})$:
\begin{eqnarray*}
I'_{\varepsilon}(u)[v] &=& \partial_{u}\mathcal I_{\varepsilon}(u, \Phi(u))[v] + \partial_{\phi}\mathcal I_{\varepsilon}(u,\Phi(u)\circ \Phi'(u)[v]\\
&=&\partial_{u}\mathcal I_{\varepsilon}(u,\Phi(u))[v]\\
&=& \int{\nabla u \nabla v} + \int{V(\varepsilon x)uv} +\int{\phi_{u}uv} - \int{f(u)u}.
\end{eqnarray*}
Then it is easy to see that  the following statements are equivalents:
\begin{itemize}
\item[i)] the pair $(u,\phi) \in H^{1}(\mathbb{R}^{3})\times \mathcal{D}$ is a critical point of $\mathcal I_{\varepsilon}$, i.e. $(u,\phi)$ is a solution of \eqref{eq:Pe*};
\item[ii)] $u$ is a critical point of $I_{\varepsilon}$ and $\phi = \phi_{u}$.
\end{itemize}

Then, solving  \eqref{eq:Pe*}  is equivalent to find critical points of $I_{\varepsilon}$, i.e., to solve 
$$-\Delta u + V(\varepsilon x)u + \phi_{u}u = f(u) \qquad \text{in} \quad \mathbb{R}^{3}.$$

Let us define de Hilbert space
$$W_{\varepsilon} = \left\{ u\in H^{1}(\mathbb{R}^{3}): \int{V(\varepsilon x)u^{2}} < +\infty \right\},$$
endowed with de scalar product and (squared) norm given by
$$(u,v)_{\varepsilon} = \int{\nabla u \nabla v} + \int{V(\varepsilon x)uv},$$
and
$$\|u\|^{2}_{\varepsilon} = \int{|\nabla u|^{2}} + \int{V(\varepsilon x)u^{2}}.$$
We will find the critical points of $I_{\varepsilon}$  in $W_{\varepsilon}$.

Defining  the Nehari manifold associated to $I_{\varepsilon}$,
$$\mathcal{N}_{\varepsilon} = \left\{ u\in W_{\varepsilon} \setminus \{0\} : J_{\varepsilon} (u)= 0\right\},$$
where
$$J_{\varepsilon}(u) = I'_{\varepsilon}(u)[u] = \|u\|^{2}_{\varepsilon} + \int{\phi_{u}u^{2}} - \int f(u)u,$$
we have the following lemma.
\begin{lemma}\label{lema.dif}
For every $u \in \mathcal{N}_{\varepsilon}, J'_{\varepsilon}(u)[u] <0$ and there are positive constants $h_{\varepsilon}, k_{\varepsilon}$, such that $\|u\|_{\varepsilon} \geq h_{\varepsilon}, I_{\varepsilon}(u) \geq k_{\varepsilon}$. Moreover, $\mathcal{N}_{\varepsilon}$ is diffeomorphic to the set $$\mathcal{S}_{\varepsilon} = \{ u\in W_{\varepsilon}: \|u\|_{\varepsilon} = 1, u> 0 \quad a.e.\}.$$
\end{lemma}
\begin{proof}
The proof follows the same steps of \cite{z4}.
\end{proof}

 By the assumptions on $f$, the functional $I_{\varepsilon}$ has the Mountain Pass geometry shown below:
\begin{itemize}
\item[(MP1)]$I_{\varepsilon}(0) = 0;$
\item[(MP2)] due to \eqref{f2} and \eqref{f3}, for all $\xi >0$ there exists $M_{\xi} >0$ such that $$F(u) \leq \xi u^{2} + M_{\xi}|u|^{q_{0}+1}.$$
Knowing that $\phi_{u}>0$ (for $u\neq0$)
\begin{eqnarray*}
I_{\varepsilon}(u)&\geq& \frac{1}{2}\|u\|^{2}_{\varepsilon} - \int{F(u)} \geq \frac{1}{2}\|u\|^{2}_{\varepsilon} - \xi \|u\|^{2}_{2}- M_{\xi}\|u\|^{q_{0}+1}_{q_{0}+1}  \\
&\geq& \frac{1}{2}\|u\|^{2}_{\varepsilon} -\xi C_{1}\|u\|^{2}_{\varepsilon} - M_{\xi} C_{2}\|u\|^{q_{0}+1}_{\varepsilon},
\end{eqnarray*}
and then, for $\|u\|^{2}_{\varepsilon} = \rho$ small enough, we  conclude that $I_{\varepsilon}$ has a strict local minimum at $u=0$.
\item[(MP3)] By \eqref{f4} we have $F(t) \geq Ct^{K}$ where  $ C>0$ and $K>4.$
Fixed $v\in C^{\infty}_{0}(\mathbb{R}^{3}), v>0$, we have $\phi_{tv} = t^{2}\phi_{v}$ 
and then
$$
I_{\varepsilon}(tv) = \frac{t^{2}}{2}\|v\|^{2}_{\varepsilon} + \frac{t^{4}}{4}\int \phi_{v}v^{2} - \int F(tv)
\leq\frac{t^{2}}{2}\|v\|^{2}_{\varepsilon} + \frac{t^{4}}{4}\int{\phi_{v}v^{2}} - Ct^{K}\int v^{K}.
$$
So, with $t$ big enough, we get $I_{\varepsilon} (tv)< 0$.
\end{itemize}

Denoting with
$$c_{\varepsilon} = \inf_{\gamma \in \mathcal{H}_{\varepsilon}} \sup_{t \in [0,1]} I_{\varepsilon}(\gamma(t)), \quad \mathcal{H}_{\varepsilon}= \{ \gamma \in C([0,1],W_{\varepsilon}):\gamma(0) = 0, I_{\varepsilon}(\gamma(1))<0 \},$$
the Mountain Pass level, and  with
$$m_{\varepsilon} = \inf_{u\in \mathcal{N}_{\varepsilon}} I_{\varepsilon}(u)$$
the ground state level, we know by \cite{willem1997minimax} that 
$$c_{\varepsilon} = m_{\varepsilon} = \inf_{u \in W_{\varepsilon} \setminus \{ 0 \} } \sup_{t\geq 0} I_{\varepsilon}(tu).$$

\subsection{The problem at infinity}
Let us consider the ``limit'' problem (the autonomous problem) associated to \eqref{eq:Pe*}, that is,
\begin{equation}\label{eq:Amu}\tag{$A_{\mu}$}
-\Delta u + \mu u = f(u) \quad \text{in } \mathbb R^{3}
\end{equation} 
where $\mu>0$ is a constant. The solutions are critical points of the functional
$$E_{\mu}(u) = \frac{1}{2}\int{|\nabla u|^{2}}+\frac{\mu}{2}\int{u^{2}} - \int{F(u)},$$
in $H^{1}(\mathbb{R}^{3}).$ We will denote with $H^{1}_{\mu}(\mathbb{R}^{3})$ simply the space $H^{1}(\mathbb{R}^{3})$ endowed with the (equivalent squared) norm
$$\|u\|^{2}_{H^{1}_{\mu}(\mathbb{R}^{3})} := \|\nabla u\|^{2}_{2} + \mu \|u\|^{2}_{2}.$$
By the assumptions of the nonlinearity $f$, it is easy to see that the functional $E_{\mu}$ has the Mountain Pass geometry (similarly to $I_{\varepsilon}$), with Mountain Pass level
$$c^{\infty}_{\mu}:= \inf_{\gamma \in \mathcal{H}_{\mu}} \sup_{t \in [0,1]} E_{\mu}(\gamma (t)),$$
$$ \mathcal{H}_{\mu}:=\left\{ \gamma \in C([0,1], H^{1}_{\mu}(\mathbb{R}^{3})): \gamma(0) = 0 , E_{\mu}(\gamma(1))<0\right\}.$$
Introducing the set
$$\mathcal{M}_{\mu}:=\left\{ u\in H^{1}(\mathbb{R}^{3})\setminus \{0\} : \|u\|^{2}_{H^{1}_{\mu}} = \int f(u)u\right\},$$
it is standard to see that (like in Lemma \ref{lema.dif}):
\begin{itemize}
	\item $\mathcal{M}_{\mu}$ has the structure of a differentiable manifold (said the Nehari manifold associated to $E_{\mu}$);
	\item $\mathcal{M}_{\mu}$ is bounded away from zero and radially homeomorfic to the subset of positive functions on the unit sphere (a kind of $\mathcal{S}_{\varepsilon}$, see Lemma \ref{lema.dif});
	\item the Mountain Pass level $c^{\infty}_{\mu}$ coincides with the ground state level 
	$$m^{\infty}_{\mu}:= \inf_{u \in \mathcal{M}_{\mu}} E_{\mu}(u) > 0.$$
\end{itemize}

In the next sections, we will mainly deal with $\mu = V_{0}$ and $\mu=V_{\infty}$, when finite.
 It is easy to see that $m_{\varepsilon} \geq m^{\infty}_{V_{0}}$.

\section{Compactness properties for $I_{\varepsilon}, E_{\mu}$ and the existence of a ground state solution}

Let us start by showing the boundedness of the Palais-Smale sequences for $E_{\mu}$ in $H^{1}_{\mu}(\mathbb{R}^{3})$ and $I_{\varepsilon}$ in $W_{\varepsilon}$. Let $\{ u_{n}\} \subset H^{1}_{\mu}(\mathbb{R}^{3})$ be a Palais-Smale sequence for $E_{\mu}$, that is, $|E_{\mu}(u_{n})| \leq C$ and $E'_{\mu}(u_{n}) \rightarrow 0$. Then, for large $n$,
\begin{eqnarray*}
E_{\mu}(u_{n}) - \frac{1}{K}E'_{\mu}(u_{n})[u_{n}] &=&  \frac{1}{2}\|u_{n}\|^{2}_{H^{1}_{\mu}} - \int F(u_{n}) - \frac{1}{K} \|u_{n}\|^{2}_{H^{1}_{\mu}} + \frac{1}{K}\int f(u_{n})u_{n}\\
&=& \left( \frac{1}{2} - \frac{1}{K} \right) \|u_{n}\|^{2}_{H^{1}_{\mu}} +\frac{1}{K}\int [f(u_{n})u_{n} -KF(u_{n})]\\
&\geq&\left( \frac{1}{2} - \frac{1}{K} \right)\|u_{n}\|^{2}_{H^{1}_{\mu}}.
\end{eqnarray*}
Since, on the other hand
\begin{eqnarray*}
\left|E_{\mu}(u_{n}) - \frac{1}{K}E'_{\mu}(u_{n})[u_{n}]\right| &\leq& |E_{\mu}(u_{n})| + \frac{1}{K}|E'_{\mu}(u_{n})|
 \|u_{n}\|_{H^{1}_{\mu}}\\
&<& C +  \frac{1}{K}|E'_{\mu}(u_{n})|\|u_{n}\|_{H^{1}_{\mu}},
\end{eqnarray*}
we conclude that   $\{ u_{n}\}$ is bounded.

Arguing similarly  
we conclude that any Palais-Smale sequence $\{u_{n}\}$  for  $I_{\varepsilon}$  is bounded in $W_{\varepsilon}$.

In order to prove compactness, we need some preliminaries lemmas. 
\begin{lemma}
If $\{u_{n}\}$ is bounded in $H^{1}(\mathbb{R}^{3})$ and for some $R>0$ and $2\leq r \leq 2^{*}=6$, we have
$$\sup_{x\in \mathbb{R}^{3}}\int_{B_{R}(x)}|u_{n}|^{r} \rightarrow 0 \quad \textrm{as} \quad n\rightarrow \infty,$$
then $u_{n} \rightarrow 0$ in $L^{p}(\mathbb{R}^{3})$ for $2<p<2^{*}$.
\end{lemma}
\begin{proof}
See \cite[Lemma I.1]{Lions}.
\end{proof}

The next results are proved as in \cite{2}.

\begin{lemma} 
Let $\{ u_{n} \} \subset W_{\varepsilon}$ be bounded and such that $I'_{\varepsilon}(u_{n}) \rightarrow 0$. Then, we have either 
\begin{itemize}
\item[a)] $u_{n} \rightarrow 0$ in $W_{\varepsilon}$, or
\item[b)] there exist a sequence $\{y_{n}\} \subset \mathbb{R}^{3}$ and constants $R, c >0$ such that $$\liminf_{n \rightarrow +\infty}\int_{B_{R}(y_{n})}{u^{2}_{n}} \geq c > 0.$$ 
\end{itemize}
\end{lemma}

In the rest of this paper, we assume, without loss of generality, that $0 \in M$, that is, $V(0) = V_{0}.$

\begin{lemma}
Assume that $V_{\infty} < \infty$ and let $\{v_{n}\} \subset W_{\varepsilon}$ be a $(PS)_{d}$ sequence for $I_{\varepsilon}$ such that $v_{n} \rightharpoonup 0$ in $W_{\varepsilon}$. Then
$v_{n} \nrightarrow 0$ in $W_{\varepsilon} $ implies $ d\geq m^{\infty}_{V_{\infty}}.$
\end{lemma}
Then the Palais-Smale condition holds:
\begin{proposition}
The functional $I_{\varepsilon}$ in $W_{\varepsilon}$ satisfies the $(PS)_{c}$ condition
\begin{itemize}
\item[1.] at any level $c<m^{\infty}_{V_{\infty}}$, if $V_{\infty}<\infty,$
\item[2.] at any level $c \in \mathbb{R}$, if $V_{\infty} = \infty.$
\end{itemize}
\end{proposition}
\begin{proof}
The proof follows from the properties of the operator $$A:u\mapsto \int{\phi_{u}u^{2}},$$
and the ideas contained in \cite{2}.
\end{proof}

Then we have

\begin{proposition}
The functional $I_{\varepsilon}$ restricted to $\mathcal{N}_{\varepsilon}$ satisfies the $(PS)_{c}$ condition:
\begin{itemize}
\item[1.] at any level $c<m^{\infty}_{V_{\infty}}$, if $V_{\infty}<\infty,$
\item[2.] at any level $c \in \mathbb{R}$, if $V_{\infty} = \infty.$
\end{itemize}
Moreover the constrained critical points of the functional $I_{\varepsilon}$ on $\mathcal{N}_{\varepsilon}$ are critical points of $I_{\varepsilon}$ in $W_{\varepsilon}$, hence solution of \eqref{eq:Pe*}.
\end{proposition}

In order to prove our main result, we recall the lemma contained in \cite{1} about the problem 
\eqref{eq:Amu}:
\begin{lemma}[Ground state for the autonomous problem]
Let $\{u_{n}\} \subset \mathcal{M}_{\mu}$ be a sequence satisfying $E_{\mu}(u_{n}) \rightarrow \mu^{\infty}_{\mu}$. Then, up to subsequences the following alternative holds:
\begin{itemize}
\item[a)] $\{u_{n}\}$ strongly converges in $H^{1}(\mathbb{R}^{3})$;
\item[b)] there exists a sequence $\{\widetilde{y}_{n}\} \subset \mathbb{R}^{3}$ such that $u_{n}(.+ \widetilde{y}_{n})$ strongly converges in $H^{1}
(\mathbb{R}^{3})$.
\end{itemize}
In particular, there exists a minimizer $\mathfrak{m}_{\mu} \geq 0$ for $m^{\infty}_{\mu}$.
\end{lemma}

Now we can prove the existence of a ground state solution for our problem. 
This is a result like   \cite[Theorem 1]{1}.
\begin{theorem}
Suppose that $V$ and $f$ verify \eqref{V} and \eqref{f1}-\eqref{f5}. Then there exists a ground state solution $\mathfrak{u}_{\varepsilon} \in W_{\varepsilon}$ of \eqref{eq:Pe*}:
\begin{itemize}
\item[1.] for every $\varepsilon \in (0,\bar{\varepsilon}]$, for some $\bar{\varepsilon}>0$, if $V_{\infty} < \infty$;
\item[2.] for every $\varepsilon > 0$, if $V_{\infty} = \infty$.
\end{itemize}
\end{theorem}
\begin{proof}
The proof follows the same lines of \cite{2} in the both cases, $V_{\infty} < \infty$ and $V_{\infty} = \infty$.
\end{proof}

\section{Proof of Theorem \ref{th:main}}
We follow the steps as in \cite{1}, to which we refer for the proofs.
Let us start with a fundamental result.
\begin{lemma}
Let $\varepsilon_{n} \rightarrow 0^{+}$ and $u_{n} \in \mathcal{N}_{\varepsilon_{n}}$ be such that $I_{\varepsilon_{n}}(u_{n}) \rightarrow m^{\infty}_{V_{\infty}}$. Then there exists a sequence $\{\widetilde{y}_{n}\} \subset \mathbb{R}$ such that $u_{n}(.+\widetilde{y}_{n})$ has a convergent subsequence in $H^{1}(\mathbb{R}^{3})$. Moreover, up to a subsequence, $y_{n}:= \varepsilon_{n}\widetilde{y}_{n} \rightarrow y \in M.$
\end{lemma}

To define the barycenter map, we first define for $\delta >0$ (later on it will be fixed conveniently), 
a smooth nonincreasing cut-off function $\eta$ in $C_{0}^{\infty}(\mathbb R^{3}, [0,1])$ such that
$$\eta(s) =\left\{
\begin{array}{rc}
1, &\mbox{if} \quad 0\le s \le {\delta}/{2},\medskip\\
0, &\mbox{if} \quad s \geq \delta.
\end{array}
\right.
$$
Let $\mathfrak{m}_{V_{0}}$ be a ground state solution of  problem
\eqref{eq:Amu}  with $\mu = V_{0}$.
For any $y\in M$, let us define
$$\Psi_{\varepsilon,y}(x):= \eta(|\varepsilon x-y|)\mathfrak{m}_{V_{0}}\left(\frac{\varepsilon x-y}{\varepsilon}\right).$$
Now, let $t_{\varepsilon}>0$ verifying $\max_{t\geq 0} I_{\varepsilon}(t\Psi_{\varepsilon,y}) = I_{\varepsilon}(t_{\varepsilon}\Psi_{\varepsilon,y})$, so that $t_{\varepsilon}\Psi_{\varepsilon,y} \in \mathcal{N}_{\varepsilon}$, and 
 define the map $\Phi_{\varepsilon}:y \in M \mapsto t_{\varepsilon}\Psi_{\varepsilon,y} \in \mathcal{N}_{\varepsilon}$.

By construction, $\Phi_{\varepsilon}(y)$ has compact support for any $y\in M$ and $\Phi_{\varepsilon}(y)$  is continuous.
The next result will help us to define a map from $M$ to a suitable sublevel in the Nehari manifold.

\begin{lemma}
The function $\Phi_{\varepsilon}$ satisfies
$$\lim_{\varepsilon \rightarrow 0^{+}} I_{\varepsilon}(\Phi_{\varepsilon}(y)) = m^{\infty}_{V_{0}},$$
uniformly in $y \in M$.
\end{lemma}
By the previous lemma, $h(\varepsilon) := |I_{\varepsilon}(\Phi_{\varepsilon}(y)) - m^{\infty}_{V_{0}}| = o(1)$ for $
\varepsilon \rightarrow 0^{+}$ uniformly in $y$, and then $I_{\varepsilon}(\Phi_{\varepsilon}(y)) - m^{\infty}_{V_{0}} \leq 
h(\varepsilon).$ In particular, the sublevel set in the Nehari 
$$\mathcal{N}^{m^{\infty}_{V_{0}} + h(\varepsilon)}
_{\varepsilon}:=\left\{u \in \mathcal{N}_{\varepsilon}: I_{\varepsilon}(u) \leq m^{\infty}_{V_{0}} + h(\varepsilon) 
\right\}$$ is not 
empty, since for sufficiently small $\varepsilon$, 
\begin{equation}
\forall y \in M: \Phi_{\varepsilon}(y) \in \mathcal{N}^{m^{\infty}_{V_{0}} + h(\varepsilon)}_{\varepsilon}.
\label{4.12}
\end{equation}
Now, we fix the $\delta >0$ mentioned before such that $M$ and $$M_{2\delta}:=\{x \in \mathbb{R}^{3}:d(x,M) \leq 2\delta\}$$
 are homotopically equivalent.

Take $\rho = \rho(\delta) >0$ such that $M_{2\delta} \subset B_{\rho}$ and define $\chi:\mathbb{R}^{3} \rightarrow \mathbb{R}^{3}$ as follows

$$\chi(x) =\left\{
\begin{array}{rc}
x, &\mbox{if}\quad |x|\le \rho,\\
\rho \frac{x}{|x|}, &\mbox{if}\quad |x| \geq \rho.\end{array}
\right.
$$
The barycenter map $\beta_{\varepsilon}$  is defined as
$$\beta_{\varepsilon}(u):= \frac{\displaystyle\int{\chi(\varepsilon x)u^{2}(x)}}{\displaystyle\int{u^{2}}} \in \mathbb{R}^{3},$$
for all $u\in W_{\varepsilon}$ with compact support. Some technical lemmas are stated now. For the proofs see e.g. \cite{1}.

\begin{lemma}
The function $\beta_{\varepsilon}$ satisfies $$\lim_{\varepsilon \rightarrow 0^{+}} \beta_{\varepsilon}(\Phi_{\varepsilon}(y)) = y$$ uniformly in $y\in M$.
\end{lemma}

\begin{lemma}
We have 
$$\lim_{\varepsilon \rightarrow 0^{+}} \sup_{u\in \mathcal{N}^{m^{\infty}_{V_{0}} + h(\varepsilon)}_{\varepsilon}} \inf_{y\in M} |\beta_{\varepsilon}(u) - y| =0.$$
\end{lemma}

Then  the proof of our main result   can be finished.
In virtue of the  above lemmas, there exist $\varepsilon^{*}>0$ such that
$$\forall \varepsilon \in (0,\varepsilon^{*}]: \sup_{u\in \mathcal{N}^{m^{\infty}_{V_{0}} + h(\varepsilon)}_{\varepsilon}} d(\beta_{\varepsilon}(u), M_{\delta})< \frac{\delta}{2}.$$ 
Let $M^{+}:=\{x\in \mathbb{R}^{3}:d(x,M) \leq 3\delta /2\}$ then homotopically equivalent to  $M$.
Now, reducing $\varepsilon^{*}>0$ if necessary, we can assume that the above lemmas and (\ref{4.12}) hold. So the composed map
$$M \xrightarrow{\Phi_{\varepsilon}} \mathcal{N}^{m^{\infty}_{V_{0}} + h(\varepsilon)}_{\varepsilon} \xrightarrow{\beta_{\varepsilon}} M^{+}$$ is homotopic to the inclusion map. 

In the case $V_{\infty} < \infty$ we eventually reduce $\varepsilon^{*}$ in such a way that also the Palais-Smale condition is satisfied in the interval $(m^{\infty}_{V_{0}}, m^{\infty}_{V_{0}} + h(\varepsilon))$.

From the properties of the Ljusternick-Schnirelamnn category we have $$cat(\mathcal{N}^{m^{\infty}_{V_{0}} + h(\varepsilon)}_{\varepsilon})\geq cat_{M^{+}}(M).$$
Thus, using Ljusternik-Schnirelmann theory, one is able to guarantee that there exists at least $\textrm{cat}_{M^{+}}(M) = \textrm{cat}(M)$ critical points of $I_{\varepsilon}$ constrained in $\mathcal{N}_{\varepsilon}$.  These critical points are solutions of our problem.

If $\textrm{cat}M >1$,  the existence of another critical point of $I_{\varepsilon}$ in $\mathcal{N}_{\varepsilon}$ follows from
the  ideas used in \cite{11}. 
The strategy is to exhibit a subset $\mathcal{A} \subset \mathcal{N}_{\varepsilon}$ such that

\begin{enumerate}
			\item $\mathcal{A}$ is not contractible in $\mathcal{N}^{m^{\infty}_{V_{0}} + h(\varepsilon)}_{\varepsilon}$ \smallskip,
			\item $\mathcal{A}$ is contractible $\mathcal{N}^{\bar{c}}_{\varepsilon} = \{u \in \mathcal{N}_{\varepsilon} :I_{\varepsilon}(u) \leq \bar{c} \}$ for some $\bar{c} > m^{\infty}_{V_{0}} + h(\varepsilon).$
			\end{enumerate}
This would implies, since the Palais-Smale holds,  the existence of a critical level between $m^{\infty}_{V_{0}} + 
h(\varepsilon)$ and $\bar{c}$.

Take $\mathcal{A} := \Phi_{\varepsilon}(M)$ which is not contractible in $\mathcal{N}^{m^{\infty}_{V_{0}} + h(\varepsilon)}_{\varepsilon}$.
Let $t_{\varepsilon}(u) >0$ the unique positive number such that $t_{\varepsilon}(u)u \in \mathcal{N}_{\varepsilon}$.

Choosing a function $u^{*} \in W_{\varepsilon}$ such that $u^{*} \geq 0, I_{\varepsilon}(t_{\varepsilon}(u^{*})u^{*}) > m^{\infty}_{V_{0}} + h(\varepsilon)$ and considering the compact and contractible cone  $$\mathfrak{C}:=\{tu^{*} + (1-t)u: t\in [0,1], u\in \mathcal{A}\},$$
we observe that, since the functions in $\mathfrak{C}$ have to be positive on a set of nonzero measure, it has to be $0\notin\mathfrak{C}$.
Then 
let $t_{\varepsilon}(\mathfrak{C}) = \{t_{\varepsilon}(w)w: w\in \mathfrak{C}\} \subset \mathcal{N}_{\varepsilon}$ and
 $$\bar{c} := \max_{t_{\varepsilon}(\mathfrak{C})} I_{\varepsilon}> m^{\infty}_{V_{0}} + h(\varepsilon).$$
It follows that $\mathcal{A} \subset t_{\varepsilon}(\mathfrak{C}) \subset \mathcal{N}_{\varepsilon}$ e $t_{\varepsilon}(\mathfrak{C})$ is contractible in $\mathcal{N}^{\bar{c}}_{\varepsilon}$.
Then there is a critical level for $I_{\varepsilon}$ greater than $m^{\infty}_{V_{0}} + h(\varepsilon)$, hence different from the previous one.

\subsection*{Acknowledgements}
The authors declare that they have no conflict of interest and contributed equally.  Bruno Mascaro was  supported by  Capes, Brazil, and G. Siciliano  was partially supported by Fapesp grant 2019/27491-0, CNPq grant 304660/2018-3, FAPDF, CAPES (Brazil) and INdAM (Italy).

\EditInfo{October 14, 2021}{November 19, 2021}{Serena Dipierro}

\end{paper}